\documentclass[12pt]{article}

\usepackage{latexsym}
\usepackage{amssymb}

\begin{document}

\title{\Large{\bf Equitable Coloring of Interval Graphs\\
and Products of Graphs}\thanks{This paper was originally
contributed to a hitherto unpublished Festschrift in honor 
of Man Keung Siu in February 2004.}}

\author{Bor-Liang Chen\\
\normalsize     Department of Business Administration\\
\normalsize     National Taichung Institute of Technology\\
\normalsize     Taichung, Taiwan 404    \\
\normalsize     {\tt E-mail:~blchen@mail.ntit.edu.tw}\\
\and
        Ko-Wei Lih\\
\normalsize     Institute of Mathematics\\
 \normalsize    Academia Sinica\\
 \normalsize    Nankang, Taipei, Taiwan 115\\
\normalsize     {\tt E-mail:~makwlih@sinica.edu.tw}\\
\and
        Jing-Ho Yan\\
\normalsize     Department of Mathematics\\
\normalsize     Aletheia University\\
\normalsize     Tamsui, Taipei, Taiwan 251\\
\normalsize     {\tt E-mail:~jhyan@email.au.edu.tw}}

\bigskip

\date{}

\maketitle

\newtheorem{theorem}{Theorem}
\newtheorem{problem}{Problem}
\newtheorem{conjecture}{Conjecture}
\newtheorem{definition}[theorem]{Definition}
\newtheorem{lemma}[theorem]{Lemma}
\newtheorem{corollary}[theorem]{Corollary}
\newcommand{\qed}{\hfill $\Box$ }
\newcommand{\proof}{\noindent{\bf Proof.}\ \ }

\baselineskip=16pt
\parindent=0.8cm

\begin{abstract}
We confirm the equitable $\Delta$-coloring conjecture for interval graphs
and establish the monotonicity of equitable colorability for them.
We further obtain results on equitable colorability about square
(or Cartesian) and cross (or direct) products of graphs.
\end{abstract}



\section{Introduction}

All graphs $G=(V,E)$
considered in this paper are finite, loopless, and without multiple edges.
Let $C_n$ and $K_n$ denote the cycle and the complete graph
on $n$ vertices, respectively. We also use $K_{x,y}$ (or $K_{x,y,z}$) to
denote the complete bipartite (or tripartite) graph with parts of sizes
$x$ and $y$ (or $x$, $y$, and $z$).
A graph $G$ is said to be {\em $k$-colorable} if there is a function
$c: V(G) \rightarrow [k]=\{0, 1, \ldots , k-1\}$ such that adjacent vertices are mapped to
distinct numbers. The function $c $ is called a proper {\em $k$-coloring} of $G$. All pre-images
of a fixed number form a so-called {\em color class}. Each color class is an independent
set, i.e., no two vertices in the same color class are adjacent.
The smallest number $k$ such that $G$
is $k$-colorable is called the {\em chromatic number} of $G$,
denoted by $\chi (G)$. A graph $G$ is said to be {\em equitably
$k$-colorable} if there is a proper $k$-coloring whose color classes
$V_{0}, V_{1}, \ldots , V_{k-1}$ satisfy the
condition $| \,|\, V_{i}| - |\, V_{j}| \,| \leqslant 1$
for all $i, j \in [k]$. The smallest integer $n$ for which $G$ is equitably
$n$-colorable is called the {\em equitable chromatic number} of $G$,
denoted by $\chi _{_{=}}(G)$. This notion of equitable colorability was
first introduced in Meyer \cite{meyer}. It is evident that $\chi (G) \leqslant \chi _{_{=}}(G)$.

Hajnal and Szemer\'{e}di \cite{hs}
shows that a graph $G$ is equitably $k$-colorable if $k \geqslant \Delta(G)+1$,
where $\Delta(G)$ denotes the maximum degree of $G$. So we may define
the parameter $\chi_{_{=}}^*(G)$ of $G$, called the {\em equitable chromatic threshold},
to be the smallest integer $n$ such that $G$ is equitable $k$-colorable
for all $k \geqslant n$.  Thus $\chi_{_{=}}^*(G) \leqslant \Delta(G)+1$.
It is obvious that $k$-colorability is monotone in the sense that $G$ is
$k$-colorable once $k \geqslant \chi(G)$.  However, equitable $k$-colorability may
fail to be monotone.  The complete bipartite graph $K_{2m+1,2m+1}$ provides
an example showing $\chi(K_{2m+1,2m+1})=\chi _{_{=}}(K_{2m+1,2m+1})=2
< \chi_{_{=}}^{*}(K_{2m+1,2m+1})=2m+2$. The following conjecture
proposed by Chen, Lih, and Wu \cite{clw} still remains open.

\bigskip

{\bf The equitable $\Delta$-coloring conjecture.}\  \  Let $G$ be a connected
graph. If $G$ is not a complete graph, or an odd cycle, or a complete bipartite
graph $K_{2m+1,2m+1}$, then $G$ is equitably $\Delta(G)$-colorable.

\bigskip

We refer the reader to a survey on equitable colorability by Lih \cite{l} for
relevant concepts and results.  The present paper supplies proofs of some statements
announced in \cite{l}. In section 2, we will confirm the equitable
$\Delta$-coloring conjecture for interval graphs and establish
the monotonicity of equitable colorability for them.  Sections 3 and 4 will handle equitable colorability of two types of graph products, namely, the square and the cross products.

\section{Interval graphs}

A graph $G(V,E)$ is called an {\em interval} graph if there exists a family
$\{I_v \mid v \in V(G)\}$ of intervals on the real line such that
$u$ and $v$ are adjacent vertices if and only if $I_u \cap I_v \ne \emptyset$.
Such a family $\{I_v \mid v \in V(G)\}$ is commonly referred to as an
{\em interval representation} of $G$. Instead of intervals of real
numbers, these intervals may be replaced by finite intervals on a linearly
ordered set.

A {\em clique} of a graph $G$ is a complete subgraph $Q$ of $G$ such that
no complete subgraph of $G$ contains $Q$ as a proper subgraph. For an interval
graph $G$, Gillmore and Hoffman \cite{gh} shows that its cliques
can be linearly ordered as $Q_0 < Q_1 < \cdots < Q_m$ so that
for every vertex $v$ of $G$ the cliques containing $v$ occur consecutively. We assign the finite interval $I_v=[Q_i,Q_j]$ in this linear order to the vertex $v$ if all the cliques containing $v$ are precisely $Q_i, Q_{i+1}, \ldots , Q_j$. Again $u$ and $v$ are adjacent if and only if $I_u \cap I_v \ne \emptyset$. We call this representation of $G$ a {\em clique path} representation of $G$. Conversely, the existence of a clique path representation implies that the graph is an interval graph.

Once a clique path representation is given, we let
{\sf left}$(v)$ and {\sf right}$(v)$ stand for the left and right
endpoint, respectively, of the interval $I_v$.  Then the following linear
order on the vertices of $G$ can be defined.  We let $u < v$ if
({\sf left}$(u)$ $<$ {\sf left}$(v)$) or ({\sf left}$(u)$ =
{\sf left}$(v)$ and {\sf right}$(u)$ $<$ {\sf right}$(v)$).
If $u$ and $v$ have the same left and right endpoints, we choose $u < v$
arbitrarily. For any three vertices $u$, $v$, and $w$ of $G$,
this linear order satisfies the following condition.
\begin{equation}\label{int}
\mbox{If $u < v < w$ and $uw \in E(G)$, then $uv \in E(G)$.}
\end{equation}
Olariu \cite{olariu} shows that the existence of a linear order satisfying (\ref{int}) characterizes interval graphs.

\begin{theorem}
Let $G$ be a connected interval graph on $n$ vertices.
If $G$ is not a complete graph, then $G$ is equitably $\Delta(G)$-colorable.
\end{theorem}

\proof\
From a clique path representation of $G$, we linearly order the vertices
of $G$ into $v_0 < v_1 < \cdots < v_{n-1}$ as defined above to satisfy
condition (\ref{int}).
Let $(a \bmod b)$ denote the remainder of $a$ divided
by $b$. Define $c(v_i) = (i \bmod \Delta(G))$ for all $v_i \in V(G)$.
It is evident that the range of $c$ contains $\Delta(G)$ colors and the
pre-images of any two colors have sizes differing by at most one.
Suppose that $v_i < v_j$ and $c(i) = c(j)$ for a pair of adjacent vertices
$v_i$ and $v_j$. It follows that $j = i + k \Delta(G)$ for some positive
integer $k$. Condition (\ref{int}) implies that $k \ngtr 1$
and $v_i$ is adjacent to $\Delta(G)$ vertices that are greater than $v_i$. However, the connectedness of $G$ implies that $v_i$ is adjacent to at least
one smaller vertex unless $i=0$.  Since the degree of $v_i$ is at most
$\Delta(G)$, it follows that the neighbors of $v_i=v_0$ are precisely
$v_1, v_2, \ldots , v_{\Delta(G)}$.

We claim that $G$ would be a complete graph on the vertices
$v_0, v_1, \ldots , v_{\Delta(G)}$.
Since $v_0 < v_1$, either {\sf left}$(v_0)$ $=$ {\sf left}$(v_1)$ or there should be a vertex $u$ in the clique
{\sf left}$(v_0)$ such that $v_0 < u < v_1$. However, the latter is impossible. We hence further have
{\sf right}$(v_0)$ $\leqslant$ {\sf right}$(v_1)$. This implies that $v_2$ is adjacent to $v_1$ since
$v_2$ is adjacent to $v_0$. Reasoning as before, we can show that
{\sf left}$(v_1)$ $=$ {\sf left}$(v_2)$ and {\sf right}$(v_1)$ $\leqslant$ {\sf right}$(v_2)$.  Arguing inductively in this way, all the vertices
$v_0, v_1, \ldots , v_{\Delta(G)}$ are shown to be mutually adjacent.
Since $G$ is connected and each vertex in $v_0, v_1, \ldots , v_{\Delta(G)}$
has degree $\Delta(G)$, our claim is true. However, this consequence is
contradicted by our assumptions. We conclude that $c$ is a proper coloring.
\qed

\bigskip

The above proof can be modified in a straightforward manner to establish the following.

\begin{corollary}
Let $G$ be a disconnected interval graph.
If $\omega(G)$, the largest size of a clique of $G$, is at most
$\Delta(G)$, then $G$ is equitably $\Delta(G)$-colorable.
\end{corollary}

\begin{theorem}
Let $G$ be an interval graph. Then $\chi_{_{=}}(G)=\chi_{_{=}}^*(G)$.
\end{theorem}

\proof\
Let $G$ have $n$ vertices. Suppose that $c$ is an equitable
$k$-coloring of $G$. Let $V_{0}, V_{1}, \ldots , V_{k-1}$ be the color classes of
$c$ such that $|\,V_j|= \lceil \frac{n-j}{k} \rceil$ for all $j \in [k]$.
We are going to modify $c$ to get an equitable $(k+1)$-coloring of $G$
by the following algorithm.

\medskip

\noindent {\bf Input.}\  \  The vertices of $G$ are listed from left to right
satisfying condition (\ref{int}).

\noindent {\bf Output.}\  \ The new color classes $V_{0}, V_{1}, \ldots , V_{k}$
are produced so that $|\,V_j|=\lceil \frac{n-j}{k+1} \rceil$ for all $j \in [k+1]$.

\noindent {\bf Initialization.}\  \ Let $S \leftarrow \{ m \mid \mbox{ $0 \leqslant m \leqslant k-1$ and $|\,V_m| > \lceil (n - m)/(k + 1) \rceil $} \}$,
$G_0 \leftarrow G \setminus \bigcup \{ V_m \mid \mbox{ $0 \leqslant m
\leqslant k-1$ and $m \notin S$} \}$, $V_{k}  \leftarrow  \emptyset$,
and $i \leftarrow 0$.
(The sequence $S$ records which old color classes
have not been reduced to the proper size.)

\noindent {\bf Procedure.}\  \

\noindent {\bf 1.}\  \  If $S = \emptyset$, then {\sf STOP}; else do the following.

\noindent {\bf 2.}\  \ Examine each vertex of $G_i$ from left to right.
While the color of a vertex occurs the first time in $G_i$, mark that vertex.
Let $v$ be the first vertex such that $c(v) = c(u)$ for a unique $u < v$.

\noindent {\bf 3.}\  \ Let $V_{k} \leftarrow V_{k} \cup \{u\}$ and $V_{c(v)} \leftarrow
V_{c(v)} \setminus \{u\}$.

\noindent {\bf 4.}\  \ If $|\,V_{c(v)}| = \lceil \frac{n-c(v)}{k+1} \rceil$,
then $G_{i+1} \leftarrow G_i \setminus (\{\mbox{all marked vertices}\} \cup V_{c(v)})$
and $S \leftarrow S \setminus \{c(v)\}$;
else $G_{i+1} \leftarrow G_i \setminus \{\mbox{all marked vertices}\}$.

\noindent {\bf 5.}\  \ Let $i \leftarrow i+1$ and {\sf GOTO} {\bf 1}.

\medskip

Now we want to prove that this algorithm is correct.

We claim that all the vertices brought to $V_{k}$ are non-adjacent.
Suppose on the contrary that there are adjacent vertices $x$ and $y$ in $V_{k}$.
We may let $x$ be brought to $V_{k}$ earlier than $y$.  From our procedure,
it implies that $x < y$ in the linear order of $G$.  When
$x$ was brought into  $V_{k}$, there was a vertex $z$ such that $x < z$ and they
both were in the same color class. The vertices appearing earlier than $z$ were
all excluded from further consideration by our procedure. So we must have $z < y$.
If $x$ and $y$ are neighbors, then condition (\ref{int}) implies that $x$ and $z$ are
neighbors, which is impossible.

Since the index $j$ is deleted from $S$ just as $|\,V_{j}| = \lceil \frac{n-j}{k+1} \rceil$
and since $n = \sum_{j=0}^{k }\lceil \frac{n-j}{k+1} \rceil$, our procedure stops if and only if
we have obtained $|\,V_j|=\lceil \frac{n-j}{k+1} \rceil$ for all $j \in [k+1]$.

When we start examining $G_{i}$, each old color class possesses at most $i$ marked
vertices. This is true because no two marked vertices have the same color in each round.
Suppose that $S$ is nonempty when we start examining $G_{i}$.
Then $|\,V_{j}| > \lceil \frac{n-j}{k+1} \rceil$ for all $j \in S$. Since after each
looping of our procedure the size of $V_{k}$ is increased by one,
we know that $|\,V_{k}|=i <  \lceil \frac{n-k}{k+1} \rceil \leqslant
\lceil \frac{n-j}{k+1} \rceil < |\,V_{j}|$ by our termination criterion above. Therefore $V_{j} \cap G_{i}$ contains
at least two unmarked vertices for every $j \in S$
and the execution of step 2 of our procedure can continue.
\qed

\bigskip

For a special subclass of interval graphs, the above monotonicity of
equitable coloring starts right from the chromatic number.
If an interval representation of an interval graph $G$ can be found so
that each interval is of unit length, then $G$ is called a {\em unit interval}
graph. A unit interval graph can be equivalently characterized as
a {\em claw-free} interval graph, i.e., an interval graph containing no
$K_{1,3}$ as an induced subgraph. A result of de Werra \cite{dw} implicitly
implies that every claw-free graph is equitably $k$-colorable for all
$k \geqslant \chi(G)$.

We now supply a simple algorithm for constructing an equitable $\chi(G)$-coloring
for a unit interval graph $G$.

The vertices of a unit interval graph $G$ can be linearly ordered
$v_0 < v_1 < \cdots < v_n$ such that each clique of $G$ consists of consecutive vertices (\cite{m}).
Define $c(v_i)=(i \bmod \omega(G))$ for all $v_i \in V(G)$.
It is evident that the range of $c$ contains $\omega(G)$ colors and the
pre-images of any two colors have sizes differing by at most one.
Suppose that $v_i < v_j$ and $c(v_i) = c(v_j)$ for a pair of adjacent vertices
$v_i$ and $v_j$. It follows that $j = i + k \omega(G)$ for some positive
integer $k$. This would imply that the set $\{v_i, v_{i+1}, \ldots , v_j\}$,
whose size is at least $\omega(G)+1$, is included in a clique. It follows
from this contradiction that $c$ is a proper coloring
of $G$ and $\chi_{_{=}}(G) \leqslant \omega(G)$. Since interval graphs are perfect graphs, we have $\omega(G)=\chi(G)$.

\section{Square products}

The {\em square product}, also known as the {\em Cartesian product},
of graphs $G_1(V_1,E_1)$ and $G_2(V_2,E_2)$ has vertex set $\{ (u,v) \mid  u\in V_1\
\mbox{and} \  v\in V_2\}$ such that $\{(u,x), (v,y)\}$ is an edge
if and only if  $(u=v \  \mbox{and} \  xy\in E_2)$ or $(x=y \  \mbox{and}\  uv\in E_1)\}$.
We denote the square product by $G_1 \Box G_2$.

\begin{theorem}\label{box}
If both $G_1$ and $G_2$ are equitably $k$-colorable, then $G_1\Box G_2$ is also
 equitably $k$-colorable.
\end{theorem}

\proof\
Let $U_0,U_1,\ldots ,U_{k-1}$ and $V_0, V_1,\ldots ,V_{k-1}$
denote the color classes of $G_1$ and $G_2$, respectively.
Suppose that we have $|\,U_0|=|\,U_1|=\cdots =|\,U_{a}|=\alpha$ for some $a \in [k]$
and the other color classes of $G_1$ are of size $\alpha-1$.
Similarly, suppose that $|\,V_0|=|\,V_1|=\cdots =|\,V_{b}|=\beta$ for some $b \in [k]$
and the other color classes of $G_2$ are of size $\beta-1$.

In the first stage, we are going to construct an auxiliary
Latin square $L=(a_{ij})$ of order $k$,
using the numbers in $[k]$ as entries. Let $q=\gcd(b,k)$ and $p=k/q$. So we
may write $b=mq$ for some $m$ such that $\gcd(m,p)=1$.
We use elements of $[k]$ to index the rows and columns of $L$.
The $(i,j)$-entry of $L$ is defined to be
$a_{i,j} =  (ib + \lfloor i/p \rfloor + j \bmod k)$.

Suppose that $a_{i,j}=ib + j=i'b  + j'=a_{i',j'}$
for $0 \leqslant i, i' \leqslant p-1$ and $0 \leqslant j, j' \leqslant q-1$.
Since $q$ divides both $b$ and $k$, it follows that $q$ divides $j-j'$,
and hence $j=j'$. This in turn implies that $p$ divides $(i-i')m$,
which is impossible.  So the upper left $p \times q$ corner of $L$, denoted by $L'$,
is filled up with the numbers in $[k]$, each occurring exactly once.
We observe that $a_{i,0}=ib + \lfloor i/p \rfloor = rb + s =
a_{r,s}$ if $i=sp+r$, where $0 \leqslant r \leqslant p-1$ and $0 \leqslant s \leqslant q-1$.
This means that the first column of $L$ is a concatenation of the
successive columns of $L'$, hence contains no repeated numbers.
As each row of $L$ is a cyclic exhibition of the numbers in $[k]$,
no repetitions in the first column imply that $L$ is a Latin square.

Now we divide $L$ into four subsquares $\left ( \begin{array}{cc}
A & B\\
C & D
\end{array} \right )$ so that the upper left corner $A$ is of order $a \times b$.
We observe that, if   $sp \leqslant i \leqslant sp+p-2$, then (i) $a_{sp+p-1,b-1}+1 = a_{sp,0}$;
(ii) $a_{i,b-1}+1 = a_{i+1,0}$.
The sequence obtained by concatenating the
$sp, sp+1, \ldots , sp+p-1$ rows of $A$, $0 \leqslant s \leqslant q-1$,
 is a cyclic exhibition of the numbers in
$[k]$, starting with the entry at $a_{sp,0}$ and each number
occurring exactly $m$ times. It implies that each number in $[k]$ occurs
in $A$ precisely $t$ or $t-1$ times, where $t=\lceil ab/k \rceil$.

Next we let $W_{z}= \bigcup \{U_i \Box V_j \mid a_{i,j}=z\}$
for each $z \in [k]$. It follows from the definition of a square product
that the following properties hold in $G_1 \Box G_2$.

(1)\  Every $U_i \Box V_j$ is an independent set.

(2)\  Every $(U_{i_1} \Box V_{j_1}) \cup (U_{i_2} \Box V_{j_2})$
is an independent set if $i_1 \ne i_2$ and $j_1 \ne j_2$.

\noindent
Consequently, every $W_{z}$ is an independent set in $G_1 \Box G_2$.
The numbers of occurrence of $z$ in the subsquares $A$, $B$, $C$,
and $D$ belong to two types: $(t, a-t, b-t, k-a-b+t)$ or
$(t-1, a-t+1, b-t+1, k-a-b+t-1)$. If $W_{z}$ belongs to the first type,
then $|\,W_{z}|=t\alpha\beta+(a-t)(\alpha\beta-\alpha)+(b-t)(\alpha\beta-\beta)
+(k-a-b+t)(\alpha\beta-\alpha-\beta+1)$. If $W_{z}$ belongs to the
second type, then $|\,W_{z}|=(t-1)\alpha\beta+
(a-t+1)(\alpha\beta-\alpha)+(b-t+1)(\alpha\beta-\beta)
+(k-a-b+t-1)(\alpha\beta-\alpha-\beta+1)$.  The difference between
the two sizes is precisely one. We conclude that
$W_{0}, W_{1}, \ldots , W_{k-1}$ form equitable color classes
for $G_1 \Box G_2$.
\qed

\begin{corollary}
Let $G_1$ have $n$ vertices and $G_2$ be $n$-colorable.
Then $G_1 \Box G_2$ is equitably $n$-colorable.
\end{corollary}

\proof\
Let the vertex set of $G_1$ be $\{u_0,u_1,$ $\ldots ,$ $u_{n-1}\}$.
Since $G_2$ is $n$-colorable, let $V_0,V_1,$ $\ldots ,$ $V_{n-1}$
be a set of color classes. Define $U_k=\bigcup \{\{u_i\} \times V_j
\mid j-i \equiv k \pmod{n}\}$ for $0 \leqslant k \leqslant n-1$.  Thus each
$U_k$ is an independent set in $G_1 \Box G_2$ and $|\,U_k|$
is equal to the order of $G_2$.
\qed

\begin{corollary}
Let $G=G_1 \Box G_2 \Box \cdots \Box G_n$, where each $G_i$ is a path,
a cycle, or a complete graph. Then we have $\chi(G)=\chi_{_{=}}(G)=\chi_{_{=}}^*(G)=\max\{\chi(G_i) \mid 1
\leqslant i \leqslant n\}$.
\end{corollary}

\proof\
This statement follows from Theorem \ref{box} together with the well-known
fact (\cite{s}) that $\chi(G)= \max \{\chi(G_i) \mid 1
\leqslant i \leqslant n\}$.
\qed

\begin{corollary}\label{coro7}
We have $\chi_{_{=}}^*(G_1 \Box G_2) \leqslant
\max\{\Delta(G_1)+1,\Delta(G_2)+1\}$.
\end{corollary}

\proof \
Since $\chi_{_{=}}^*(G_1) \leqslant \Delta(G_1)+1$ and $\chi_{_{=}}^*(G_2)
\leqslant \Delta(G_2)+1$, we have $\chi_{_{=}}^*(G_1 \Box G_2) \leqslant \max
\{\Delta(G_1)+1,\Delta(G_2)+1\}$ by Theorem \ref{box}.
\qed

\begin{corollary}
Suppose that $G_1$ and $G_2$ are graphs each with at least one edge.
Then $G_1 \Box G_2$ is equitably $\Delta(G_1 \Box G_2)$-colorable.
\end{corollary}

\proof\
Since neither $G_1$ nor $G_2$ consists of isolated vertices,
we have $\Delta(G_1 \Box G_2) \geqslant \max\{\Delta(G_1)+1,\Delta(G_2)+1\}$.
Corollary \ref{coro7} implies that $\chi_{_{=}}^*(G_1 \Box G_2) \leqslant
\max\{\Delta(G_1)+1,\Delta(G_2)+1\} \leqslant \Delta(G_1 \Box G_2)$.
\qed

\bigskip

If we weaken the assumption on $G_1$ in Theorem \ref{box}
to that of its $k$-colorability, then the conclusion may not follow.
Let $K_{1,5}$ denote the star graph on $6$ vertices and $P_3$
the path on $3$ vertices. The cross product $K_{1,5} \Box P_3$
is a bipartite graph with one part $A$ of size $7$ and the other part
$B$ of size $11$. Let us consider any proper $2$-coloring of this product.
Since there is a vertex $x$ in $B$ that is adjacent to every vertex in $A$,
none of the vertices in $A$ belong to the color class containing $x$.
But any vertex in $B$ is adjacent to some vertex in $A$. Therefore,
this $2$-coloring cannot be equitable.  This example shows that,
even if $\chi(G_1)=\chi_{_{=}}(G_2)=k$, the product $G_1 \Box G_2$
may not be equitable $k$-colorable.

If we assume that $\chi_{_{=}}(G_1)$ $=$ $\chi_{_{=}}(G_2)=k$, it may not
lead to the conclusion $\chi_{_{=}}(G_1 \Box G_2)$ $=k$.
Let us consider $K_{1,2n} \Box K_{1,2n}$. Let the vertex set
of $K_{1,2n}$ be $\{a_0, a_1, \ldots , a_{2n}\}$ so that $a_0$ is
the vertex of degree $2n$. It is easy to see that $\chi_{_{=}}(K_{1,2n})=n+1$.
The following array gives an equitable $4$-coloring of $K_{1,2n} \Box K_{1,2n}$.
(The entry at position $(i,j)$ is the color given to the vertex $(a_i,a_j)$.)

\medskip

\[
\begin{array}{cccccc}
  &     & 0 & \overbrace{\begin{array}{ccc} 1  &  \cdots & 1 \end{array}}^{n}
& 2 & \overbrace{\begin{array}{ccc} 2 & \cdots & 2 \end{array}}^{n-1}\\
n  &  \Bigg\{ &
     {\begin{array}{c}
        3\\
        \cdot \\
        3 \end{array}}
         & {\begin{array}{ccc}
        0 &  \cdots &  0 \\
        \cdot & \cdots & \cdot \\
        0 &  \cdots &  0  \end{array}}
          & {\begin{array}{c}
        0\\
        \cdot \\
        0 \end{array}}
          & {\begin{array}{ccc}
        1 &  \cdots &  1 \\
        \cdot & \cdots & \cdot \\
        1 &  \cdots &  1  \end{array}} \\
n  &  \Bigg\{ &
     {\begin{array}{c}
        1\\
        \cdot \\
        1 \end{array}}
         & {\begin{array}{ccc}
        2 &  \cdots &  2 \\
        \cdot & \cdots & \cdot \\
        2 &  \cdots &  2  \end{array}}
          & {\begin{array}{c}
        3\\
        \cdot \\
        3 \end{array}}
          & {\begin{array}{ccc}
        3 &  \cdots &  3 \\
        \cdot & \cdots & \cdot \\
        3 &  \cdots &  3  \end{array}}
\end{array}
\]

\bigskip

The following example shows that $\chi_{_{=}}(G_1 \Box G_2) \leqslant
\max\{\chi_{_{=}}(G_1), \chi_{_{=}}(G_2)\}$ is false in general.
Let $G_1=K_{3,3}$ and $G_2=K_{2,1,1}$.
We have $\chi_{_{=}}(G_1)=2$ and $\chi_{_{=}}(G_2)=3$,
but $\chi_{_{=}}(G_1 \Box G_2)=4$. It is easy to see that
$G_1 \Box G_2$ is equitably $4$-colorable. We want to show that it
is not equitably $3$-colorable.
We write the vertices of $G_1$ into a sequence $[u_0,u_1,u_2,v_0,v_1,v_2]$
so that $\{u_0,u_1,u_2\}$ and $\{v_0,v_1,v_2\}$ form independent sets, respectively.
We write the vertices of $G_2$ into a sequence $[a_0,a_1,b,c]$ so that,
except $a_0$ and $a_1$, all pairs of vertices are adjacent.  Now we arrange the
vertices of $G_1 \Box G_2$ into a $6 \times 4$ array. Suppose that there
were an equitable $3$-coloring of this array.  Thus every color class contains
exactly $8$ vertices.  Each pair $(x,a_0)$ and $(x,a_1)$ must have the same color
since they are adjacent to the two endpoints of the edge $(x,b)(x,c)$. It implies that the
first column has at least two colors. Since we cannot have a pair
$(u_i,a_0)$ and $(v_j,a_0)$ with the same color, either all $(u_i,a_0)$'s
are of the same color or all $(v_i,a_0)$'s are of the same color.
Either possibility implies that some color class would contain
$9$ vertices.

In general, let $G_1$ be equitably $k_1$-colorable and $G_2$ be
equitably $k_2$-colorable.  It remains open to find conditions that force
$G_1 \Box G_2$ to be equitably $(\max\{k_1,k_2\})$-colorable.

\section{Cross products}

The {\em cross product}, also known as the {\em direct product},
of graphs $G_1(V_1,E_1)$ and $G_2(V_2,E_2)$ has vertex set
$\{ (u,v) \mid  u\in V_1\  \mbox{and} \  v\in V_2\}$ such that
$\{(u,x), (v,y)\}$ is an edge if and only if  $uv\in E_1$ and
$xy\in E_2$. We denote the cross product by $G_1 \times G_2$.

\begin{lemma}\label{time}
We have $\chi_{_{=}}(G_1 \times G_2) \leqslant \min\{|\,V(G_1)|,|\,V(G_2)|\}$.
\end{lemma}

\proof\ Let $V(G_1)=\{u_0,u_1, \ldots ,u_m\}$ and
$U_i=\{u_i\} \times V(G_2)$ for all $0 \leqslant i\leqslant m$.
Then $U_i$ is an independent set of $G_1 \times G_2$ and
$|\,U_i|=|\,V(G_2)|$ for every $0 \leqslant i \leqslant m$.
Thus $\chi_{_{=}}(G_1 \times G_2) \leqslant |\,V(G_1)|$.
Similarly, we have $\chi_{_{=}}(G_1 \times G_2) \leqslant |\,V(G_2)|$.
\qed

\begin{corollary}
We have $\chi_{_{=}}(K_m \times K_n)=\min\{m,n\}$.
\end{corollary}

\proof\
Duffus, Sands, and Woodrow \cite{dsw} shows
that $\chi(K_m \times K_n)=\min\{\chi (K_m),\chi(K_n)\}$.
Then Lemma \ref{time} implies the result.
\qed

\bigskip

We note that $\chi_{_{=}}^*(G_1 \times G_2) \leqslant \min\{|\,V(G_1)|,|\,V(G_2)|\}$ is false in general.
For instance, $K_2\times K_n$ is not equitably $(\frac{n+1}{2})$-colorable if
$n > 1$ and $n \equiv 1 \pmod{4}$.

Let $n=4k+1$ for some $k \geqslant 1$. We observe that
$K_2\times K_n$ is equal to $K_{n,n} - M$ for a complete matching $M$.
If there were an equitable $(\frac{n+1}{2})$-coloring of $K_2\times K_n$, then
there would be two color classes of size $3$ and $2k-1$ color classes
of size $4$.  Any color class of size $3$ cannot contain vertices in both
parts. Yet neither $4k-2$ nor $4k-5$ is divisible by $4$. Hence the
desired equitable color classes cannot exist.

We would surmise that $\chi_{_{=}}^*(G_1 \times G_2)$ $\leqslant$ $\max\{|\,V(G_1)|,|\,V(G_2)|\}$ should be true.

\begin{theorem}
Let $m,n \geqslant 3$. Then
\[
\chi_{_{=}}(C_m \times C_n)=\chi_{_{=}}^*(C_m \times C_n)
=\left\{  \begin{array}{ll}
2, & \mbox{if  $mn$ is even;}\\
3, & \mbox{otherwise.}
\end{array} \right.
\]
\end{theorem}

\proof\
Let $C_m$ be the cycle $u_0u_1\cdots u_{m-1}u_0$ and $C_n$ be the cycle $v_0v_1\cdots v_{n-1}v_0$.
We note that $C_m \times C_n$ is a $4$-regular graph. Hence it is equitably
$k$-colorable for all $k \geqslant 5$.

{\em Case 1.} \  We use two colors.

If $mn$ is even, then $C_m \times C_n$ is a bipartite graph
with parts of equal size. Hence $\chi_{_{=}}(C_m \times C_n)=2$.
If $m \leqslant n$ are both odd, then there exists an odd cycle in $C_m \times C_n$:
$(u_0,v_0)(u_1,v_1)$ $\cdots$ $(u_{m-1},v_{m-1})$$(u_{m-2},v_{m})$
$(u_{m-1},v_{m+1})$ $\cdots$ $(u_{m-1},v_{n-1})(u_0,v_0)$.
Hence $\chi_{_{=}}(C_m \times C_n) \geqslant 3$.

{\em Case 2.} \  We use three colors.

It is straightforward to verify
the colorings to be defined in the following subcases are equitable
$3$-colorings of $C_m \times C_n$.

{\em Subcase 2.1.}\
Assume that $m$ or $n$, say $n$, is divisible by $3$.
Define the coloring $a(u_i,v_j) = (j \bmod 3)$.

{\em Subcase 2.2.}\
Assume that $m-1$ or $n-1$, say $n-1$, is divisible by $3$.
Also assume that $n > 4$. Define the coloring
\[
b(u_i,v_j)=\left\{  \begin{array}{ll}
0, & \mbox{if $j=n-2$;}\\
1, & \mbox{if $j=n-1$ or ($j=n-4$ and $i < \lceil m/3 \rceil$);}\\
2, & \mbox{if $j=n-3$ or ($j=0$ and $i < \lfloor m/3 \rfloor$);}\\
a(u_i,v_j), & \mbox{otherwise.}
\end{array} \right .
\]

{\em Subcase 2.3.}\
Assume that $m-2$ or $n-2$, say $n-2$, is divisible by $3$.
Also assume that $n > 5$. Define the coloring
\[
c(u_i,v_j)=\left\{  \begin{array}{ll}
0, & \mbox{if $j=n-3$;}\\
2, & \mbox{if  $j=n-2$, or ($j=0$ and $i < \lceil m/3 \rceil$),}\\
     & \mbox{or ($j=n-4$ and $i < \lfloor m/3\rfloor$);}\\
a(u_i,v_j), & \mbox{otherwise}.
\end{array} \right .
\]

{\em Subcase 2.4.}\
There are three remaining cases that are solved
by the following arrays of colorings.

\medskip

\[
\begin{array}{cccc}
\begin{array}{cccc}
0 & 0 & 0 & 0\\
2 & 1 & 2 & 1\\
0 & 0 & 2 & 1\\
2 & 1 & 2 &1
\end{array}

&
\qquad
\begin{array}{ccccc}
0 & 1 & 1 & 1 & 1\\
0 & 2 & 2 & 0 & 2\\
0 & 1 & 1 & 0 & 1\\
2 & 2 & 2 & 0 & 2
\end{array}

&
\qquad
\begin{array}{ccccc}
2 & 0 & 2 & 1 & 2\\
1 & 0 & 2 & 1 & 0\\
1 & 0 & 2 & 1 & 2\\
1 & 0 & 2 & 1 & 0\\
1 & 0 & 2 & 1 & 0
\end{array}
\end{array}
\]

\medskip

{\em Case 3.} \  We use four colors.

Again it is straightforward to verify
the colorings to be defined in the following subcases are equitable
$4$-colorings of $C_m \times C_n$.

{\em Subcase 3.1.}\
Assume that $m$ or $n$, say $n$, is divisible by $4$.
Define $d(u_i,v_j) = (j \bmod 4)$.

{\em Subcase 3.2.}\
Assume that $m-1$ or $n-1$, say $n-1$, is divisible by $4$.
Also assume that $n > 5$. Define the coloring
\[
e(u_i,v_j)=\left\{  \begin{array}{ll}
0, & \mbox{if $j=n-2$;}\\
1, & \mbox{if ($j=n-1$ and $i \geqslant \lfloor (m+3)/4 \rfloor$)}\\
     & \mbox{or ($j=n-5$ and $i < \lceil m/2 \rceil$);}\\
2, & \mbox{if $j=n-5$ and $i \geqslant \lceil m/2 \rceil+\lfloor m/4 \rfloor$;}\\
3,  & \mbox{if $j=n-4$ or ($j=n-1$ and $i < \lfloor (m+3)/4 \rfloor)$;}\\
d(u_i,v_j), & \mbox{otherwise.}
\end{array} \right .
\]

{\em Subcase 3.3.}\
Assume that $m-2$ or $n-2$, say $n-2$, is divisible by $4$. Define the coloring
\[
f(u_i,v_j)=\left\{  \begin{array}{ll}
1, & \mbox{if $j=n-2$ and $i \geqslant \lfloor m/2 \rfloor$;}\\
2, & \mbox{if $j=n-1$ and $i < \lfloor m/2 \rfloor$;}\\
3,  & \mbox{if $j=n-1$  and $i \geqslant \lfloor m/2 \rfloor$;}\\
d(u_i,v_j), & \mbox{otherwise.}
\end{array} \right .
\]

{\em Subcase 3.4.}\
Assume that $m-3$ or $n-3$, say $n-3$, is divisible by $4$. Define the coloring
\[
g(u_i,v_j)=\left\{  \begin{array}{ll}
1, & \mbox{if ($j=n-1$ and $i \geqslant m- \lfloor (m+2)/4 \rfloor$);}\\
     & \mbox{or ($j=n-3$ and $i \geqslant m- \lfloor (m+2)/4 \rfloor$);}\\
3,  & \mbox{if $j=n-2$  and ($i < \lceil m/4 \rceil$ or $i \geqslant \lceil m/2 \rceil$);}\\
d(u_i,v_j), & \mbox{otherwise.}
\end{array} \right .
\]

{\em Subcase 3.5.}\ There are three remaining cases that will be solved
by the following arrays of colorings.

\medskip

\[
\begin{array}{ccc}
\begin{array}{ccc}
0 & 1 & 1\\
0 & 3 & 3\\
0 & 2 & 2
\end{array}

&
\qquad
\begin{array}{ccccc}
0 & 0 & 2 & 1 & 3\\
0 & 1 & 2 & 1 & 2\\
0 & 3 & 3 & 1 & 2\\
\end{array}

&
\qquad
\begin{array}{ccccc}
0 & 1 & 2 & 0 & 1\\
0 & 1 & 2 & 3 & 1\\
0 & 1 & 2 & 3 & 1\\
0 & 3 & 2 & 3 & 2\\
0 & 3 & 2 & 3 & 2
\end{array}
\end{array}
\]

\medskip

\begin{theorem}
We have
$\chi_{_{=}}(K_n \times K_{n,n-1})=\chi_{_{=}}^*(K_n \times K_{n,n-1})=n$.
\end{theorem}

\proof\  The statement is trivial when $n=2$. Let us assume $n \geqslant 3$.
Denote the vertices of $K_n$ by $u_0,u_1, \ldots ,u_{n-1}$ and the vertices
of $K_{n,n-1}$ by disjoint parts: $A: a_0, a_1, \ldots , a_{n-1}$ and
$B: b_0,b_1, \ldots ,b_{n-2}$. We arrange vertices of $K_n \times K_{n,n-1}$
into an $n$ by $2n-1$ array so that the $i$-th row is equal to
$(u_i,a_0)(u_i,a_1)$ $\cdots$ $(u_i,a_{n-1})(u_i,b_0)(u_i,b_1)$ $\cdots$ $(u_i,b_{n-2})$.

{\em Claim 1.\
The graph $K_n\times K_{n,n-1}$ is equitably $k$-colorable for
all $k \geqslant n$.}

Let $k \geqslant n$. We are trying to equitably color $K_n\times K_{n,n-1}$
with $k$ colors.
The size of each color class should be $m$ or $m+1$, where $m=\lfloor n(2n-1)/k \rfloor \leqslant 2n-1$.
If $m=2n-1$, then the $n$ rows form an equitable $n$-coloring. Let us assume
$m \leqslant 2n-2$ and $\alpha(m+1)+\beta m=n(2n-1)$ for some $\alpha$ and
$\beta$ with $\alpha + \beta = k$.  We are going to partition the vertices into
independent sets of appropriate sizes and numbers.

We remove initial segments of length $m+1$ from successive rows
in a cyclic fashion. Once the $\alpha$ independent
sets of size $m+1$ have been removed, we partition the remaining
part of each row into segments of length $m$. After all this is done,
the number of vertices left in each row is less than $m$, hence the
second coordinates all belong to $B$.  All these leftover vertices
form an independent set. We just partition them further into
subsets of size $m$.

{\em Claim 2.\
The graph $K_n\times K_{n,n-1}$ is not equitably $k$-colorable for
any $k < n$.}

Suppose that it were equitably $k$-colorable for some $k < n$.
Then the size of each color class is at least $\lfloor n(2n-1)/k \rfloor$.
Now $n(2n-1)/k \geqslant n(2n-1)/(n-1) = 2n+1+ \frac{1}{n-1}$.
It follows that $\lfloor n(2n-1)/k \rfloor \geqslant 2n+1$.
If a color class contains two vertices whose second coordinates
belong to different parts of $K_{n,n-1}$, then their first
coordinates must equal.  However, there are at most $2n-1$
vertices with the same first coordinates. Hence the second
coordinates of a color class must come from the same part
of  $K_{n,n-1}$.

Suppose that the part having $n$ vertices is partitioned into $x$ color classes
and the part having $n-1$ vertices is partitioned into $y$ color classes.
The sizes of color classes satisfy
$|\,\lfloor n^2/x \rfloor - \lceil n(n-1)/y \rceil \,| \leqslant 1$, which
in turn implies $|\, n^2/x - n(n-1)/y \,| \leqslant 1$.
If $x \leqslant y$, then $n^2/x > (n^2-n)/x+1 \geqslant n(n-1)/y+1$.
If $x > y$, then $n > 2y$. It follows that $(n+y-1)(n-y) > n(y+1)$, and
hence $n(n-1)/y > n^2/(y+1)+1 \geqslant n^2/x +1$.
\qed

\bigskip

We note that, even if both $G_1$ and $G_2$ are equitably
$k$-colorable, $G_1 \times G_2$ may not be equitably $k$-colorable.
Let us consider $K_{m,m-1} \times K_{n,n-1}$. This is a disjoint
union of $K_{mn,(m-1)(n-1)}$ and $K_{(m-1)n,m(n-1)}$.
If we properly color this union by two colors, then every
part should be entirely colored with one color and two
parts of the same connected component should be colored with
different colors. However, the combined size of any two independent parts
is different from $2mn-m-n$. Therefore,
this disjoint union is not equitably $2$-colorable.
In particular, $\chi_{_{=}}^*(K_{3,2} \times K_{3,2}) > 2$.
(Actually, $\chi_{_{=}}^*(K_{3,2} \times K_{3,2}) = 3$.)
However, $\chi_{_{=}}^*(K_{3,2}) = 2$ shows that the inequality
$\chi_{_{=}}^*(G_1 \times G_2) \leqslant \max \{\chi_{_{=}}^*(G_1),
\chi_{_{=}}^*(G_2)\}$ is false in general.

We conclude this paper by posing the determination of the exact values for $\chi_{_{=}}^*(K_m \times K_n)$
and $\chi_{_{=}}^*(K_{m,m-1} \times K_{n,n-1})$ as an open problem.

\end{document}